\documentclass[11pt,twoside]{article}
\usepackage{a4wide}
\usepackage{amssymb}
\usepackage{amsmath}

\newtheorem{theorem}{Theorem}
\newtheorem{thkh}{ Theorem K (1926).}

\newtheorem{thga}{ Theorem G (1962).}

\newtheorem{thrynne}{Theorem R (1996).}

\newtheorem{thbd}{Theorem BD (1978).}

\newtheorem{theoremKM}{Theorem}

\newtheorem{theoremKMd}{Theorem}

\newtheorem{thsix}{Theorem \ref{t6}*}

\newtheorem{lemma}{Lemma}
\newtheorem{definition}{Definition}

\newcommand{\bp}{\mathbf{p}}
\newcommand{\f}{\mathbf{f}}
\newcommand{\bu}{{\bf u}}
\newcommand{\cRp}{\cR_{\cC}(\Phi)}
\newcommand{\rao}{R_{\alpha,1}}
\newcommand{\rak}{R_{\alpha,k}}

\newcommand{\cJp}{\mathcal{J}_{\mathcal{C}}(\Phi)}
\renewcommand{\Bbb}[1]{\mathbb{#1}}
\newcommand{\N}{{\Bbb N}}         
\newcommand{\Q}{{\Bbb Q}}         
\newcommand{\R}{{\Bbb R}}         
\newcommand{\Rp}{{\Bbb R}^{+}}    
\newcommand{\Z}{{\Bbb Z}}         
\newcommand{\ba}{\beta_\alpha}

\newcommand{\cC}{{\cal C}}

\newcommand{\cJ}{{\cal J}}

\newcommand{\cM}{{\cal M}}

\newcommand{\cQ}{{\cal Q}}
\newcommand{\cR}{{\cal R}}
\newcommand{\cS}{{\cal S}}
\newcommand{\cSM}{{\cal S}^*}
\newcommand{\cLM}{{\cal L}^*}

\newcommand{\gQ}{\mathcal{Q}}

\newcommand{\ve}{\varepsilon}

\newcommand{\p}{\psi}

\newcommand{\diam}{\operatorname{diam}}

\newcommand{\set}[1]{\left\{#1\right\}}
\newcommand{\vv}[1]{{\mathbf{#1}}}

\renewcommand{\le}{\leqslant}
\renewcommand{\ge}{\geqslant}
\newcommand{\ra}{R_{\alpha}}

\begin{document}

\title{A note on simultaneous  Diophantine approximation on  \\ planar curves}

\author{Victor Beresnevich\footnote{Research supported by EPSRC Grant R90727/01}
\\ {\small\sc Minsk} \and  Sanju Velani\footnote{Royal Society University Research
Fellow} \\ {\small\sc York}}

\date{}

\maketitle

\centerline{{\it For Iona and  Ayesha on $N^o \!\!$ 3}}

\vspace{16mm}

\begin{abstract} Let $\cS_n(\psi_1,\dots,\psi_n)$
denote the set  of simultaneously
$(\psi_1,\dots,\psi_n)$--approximable points in $\R^n$ and
$\cSM_n(\psi)$ denote the set  of multiplicatively
$\psi$--approximable points in $\R^n$. Let $\cM$ be a  manifold in
$\R^n$. The aim is to develop  a metric theory for the sets $ \cM
\cap \cS_n(\psi_1,\dots,\psi_n) $ and $ \cM \cap \cSM_n(\psi) $
analogous to the classical theory in which   $\cM$ is simply
$\R^n$. In this note, we mainly restrict our attention to the case
that $\cM$ is a planar curve $\cC$.  A complete Hausdorff
dimension theory is established for the sets $ \cC \cap
\cS_2(\psi_1,\psi_2) $ and $ \cC \cap \cSM_2(\psi) $.  A divergent
Khintchine type result is obtained for $\cC \cap
\cS_2(\psi_1,\psi_2) $; i.e. if a certain sum diverges then the
one--dimensional Lebesgue measure on $\cC$ of $\cC \cap
\cS_2(\psi_1,\psi_2) $ is full. Furthermore, in the case that
$\cC$ is a rational quadric the convergent Khintchine type result
is obtained for both types of approximation.  Our results for $\cC
\cap \cS_2(\psi_1,\psi_2) $ naturally  generalize the dimension
and Lebesgue measure statements  of \cite{BDV03}. Within the
multiplicative framework, our results for $ \cC \cap \cSM_2(\psi)
$  constitute the first of their type.
\end{abstract}

\vspace{2cm}

\bigskip


\newpage

\section{Introduction}

\subsection{Background: two types of simultaneous approximation \label{background}}

Throughout $\psi:\R^+\to\R^+$ will denote a real, positive
decreasing function and will be referred to as an
\emph{approximating function}. Given approximating functions
$\psi_1,\dots,\psi_n$, a point $\vv y=(y_1,\dots,y_n)\in\R^n$ is
called {\it simultaneously $(\psi_1,\dots,\psi_n)$--approximable} if
there are infinitely many $q\in\N$ such that $$ \|q y_i\|<\psi_i(q)
\hspace{17mm}  1\le i\le n  $$ where $\|x\|=\min\{|x-m|:m\in\Z\}$.
In the case $\psi_i : h \to h^{-v_i}$ with $v_i > 0$,  the point
$\vv y$ is said to be {\it simultaneously
$(v_1,\dots,v_n)$--approximable}.  The set of simultaneously
$(\psi_1,\dots,\psi_n)$--approximable points in $\R^n$ will be
denoted by  $\cS_n(\psi_1,\dots,\psi_n)$ and similarly
$\cS_n(v_1,\dots,v_n)$ will denote the set of simultaneously
$(v_1,\dots,v_n)$--approximable points in $\R^n$. Geometrically,
$\vv y \in \cS_n(\psi_1,\dots,\psi_n)$ if it lies in infinitely many
$n$--dimensional `rectangular' regions centred at rational points
with `size' determined by $ \psi_1,\dots,\psi_n $.

Next, given an approximating function $\psi$, a point $\vv
y=(y_1,\dots,y_n)\in\R^n$ is called {\it multiplicatively
$\psi$--approximable} if there are infinitely many $q\in\N$ such
that $$  \prod_{i=1}^n \|qy_i\| \ < \ \psi(q) \ .
  $$  In the case $\psi
: h \to h^{-v}$ with $v > 0$ the point $\vv y$ is said to be {\it
multiplicatively $v$--approximable}.  The set of multiplicatively
$\psi$--approximable points in $\R^n$ will be denoted by
$\cSM_n(\psi)$ and similarly $\cSM_n(v)$ will denote the set of
multiplicatively $v$--approximable points in $\R^n$. In this
multiplicative setup, $\vv y \in \cSM_n(\psi)$ if it lies in
infinitely many $n$--dimensional `hyperbolic' regions centred  at
rational points with `size' determined by $ \psi$.

It is readily verified that

\begin{equation}\label{s:004}
  \cS_n(\psi_1,\dots,\psi_n)\subset\cSM_n(\psi)\quad\text{whenever}\quad\psi\ge\psi_1\cdots\psi_n,
\end{equation}
\begin{equation}\label{s:005}
 \hspace{4ex}  \cS_n(v_1,\dots,v_n)\subset\cSM_n(v) \hspace{1ex} \quad\text{whenever}\quad
v\le v_1+\dots+v_n.
\end{equation}

\noindent Also, in view of Minkowski's linear forms theorem which
gives rise to a general $n$--dimensional version of Dirichlet's
theorem,
\begin{equation}\label{sv:004}
\cS_n(v_1,\dots,v_n)=\R^n \hspace{9mm}  {\rm  if } \hspace{9mm}
v_1+\dots+v_n \leq 1 \ . \end{equation}
 This together with $
(\ref{s:005})$ implies that
\begin{equation}\label{sv:005}
\cSM_n(v)=\R^n  \hspace{9mm}  {\rm  if } \hspace{9mm} v \leq 1 \ .
\end{equation}

\vspace{14pt}

\noindent{\em The Lebesgue theory. \ }
 The following key results
provide beautiful and simple criteria for the  `size' of the sets
$\cS_n(\psi_1,\dots,\psi_n)$ and $\cSM_n(\psi)$ expressed in terms
of $n$--dimensional Lebesgue measure $|\  \ |_{\R^n}$. The first
is due to Khintchine \cite{Kh} and the second is due to Gallagher
\cite{gal}.

\begin{thkh}
\label{gs} Let   $\psi$ be an approximating function. Then
$$|\cS_n(\psi_1,\dots,\psi_n)|_{\R^n} =\left\{\begin{array}{ll} \mbox{\rm
Z{\scriptsize ERO}} & {\rm if} \;\;\; \sum \;  \p_1(h)\dots \p_n(h)
\;\; <\infty\\ &
\\ \mbox{\rm F{\scriptsize ULL}} & {\rm if} \;\;\; \sum  \;  \p_1(h)\dots \p_n(h) \;\;
 =\infty \; \;
\end{array}\right..$$
\end{thkh}

This theorem is a generalization of Khintchine's 1924 result which
deals with the special case  $\psi_1= \psi_2= \dots =\psi_n$.

\begin{thga}
\label{gm} Let   $\psi$ be an approximating function. Then
$$|\cSM_n(\psi)|_{\R^n} =\left\{\begin{array}{ll} \mbox{\rm
Z{\scriptsize ERO}} & {\rm if} \;\;\; \sum \;  \p(h)^n \ (\log
h)^{n-1} \;\; < \ \infty\\ &
\\ \mbox{\rm F{\scriptsize ULL}} & {\rm if} \;\;\; \sum  \;  \p(h)^n \ (\log
h)^{n-1}  \;\;
 = \ \infty \; \;
\end{array}\right..$$
\end{thga}

\noindent Here `full' simply means that the complement of the set
under consideration is of `zero' measure. Thus the
$n$--dimensional Lebesgue measure of the sets in question  satisfy
a `zero-full' law. The divergence parts of the above statements
constitute the main substance of the theorems. The convergence
parts are a simple consequence of the Borel-Cantelli lemma from
probability theory. Trivially, the convergence parts imply that
$$|\cS_n(v_1,\dots,v_n)|_{\R^n} = 0  \hspace{9mm}  {\rm  if }
\hspace{9mm} v_1+\dots+v_n >  1 \
$$and
$$|\cSM_n(v)|_{\R^n} = 0   \hspace{9mm}  {\rm  if }
\hspace{9mm} v >  1 \ .  $$  Note that   the former statement
is in fact  a consequence of the latter and (\ref{s:005}). In the
case that the set in question is of Lebesgue  measure zero, a more
delicate attribute of the `size' of the set is its Hausdorff measure
and dimension. In this article  we shall only be concerned with the
dimension theory.  The Hausdorff dimension of a set $X \in \R^n$ is
defined as follows. For $\rho > 0$, a countable collection $
\left\{B_{i} \right\} $ of Euclidean balls  in  $\R^{n}$ with
diameter $\diam(B_i)  \leq \rho $ for each $i$  such that $X \subset
\bigcup_{i} B_{i} $ is called a $ \rho $-cover for $X$.  Let $s$ be
a non-negative number and define $
 {\cal H}^{s}_{ \rho } (X)
  \; = \; \inf \left\{ \sum_{i} \diam(B_i)^s
\ :   \{ B_{i} \}  {\rm \  is\ a\  } \rho {\rm -cover\  of\ } X
\right\} $,  where the infimum is taken over all possible $ \rho
$-covers of $X$. The {\it Hausdorff dimension} dim $X$ of $X$ is
defined by  infimum over  $s$ for which $  \sup_{ \rho
> 0} {\cal H}^{s}_{ \rho } (X)
$  is zero.

\vspace{14pt}

\noindent{\em The dimension theory. \ } The following relatively
recent results provide exact formulae for the `size' of the sets
$\cS_n(v_1,\dots, v_n)$ and $\cSM_n(v)$ expressed in terms of
Hausdorff dimension. The first is due to Rynne \cite{rynne} and
the second is due to Bovey $\&$ Dodson \cite{bd}.

\begin{thrynne}
\label{rynne} Let $ v_1\geq  v_2 \geq  \dots \geq   v_n $ and $ v_1
+ v_2 + \dots +  v_n \geq 1 $. Then
$$
\dim \cS_n(v_1,\dots,v_n) \ = \ \min_{1\leq k \leq n }  \left\{
\frac{n+1 + \textstyle{\sum_{i=k}^{n}} (v_k- v_i) }{ 1+ v_k }
\right\} \ . $$
\end{thrynne}

In the case $ v_1= v_2 = \dots =  v_n $, the above statement reduces
to the classical Jarn\'{\i}k--Besicovitch theorem.

\begin{thbd}
\label{bd} Let $v \geq 1 $. Then
$$
\dim \cSM_n(v) \ = \ n-1 \, + \,  \frac{2}{v+1} \  \ .
$$
\end{thbd}

\subsection{Simultaneous approximation restricted to manifolds \label{simonman}}

Let $\cM$ be a  manifold in $\R^n$. In short, the aim is to
develop  a metric theory for the sets $ \cM \cap
\cS_n(\psi_1,\dots,\psi_n) $ and $ \cM \cap \cSM_n(\psi) $
analogous to that described above in which $\cM$ is simply $\R^n$.
The fact that the points $\vv y$ of interest consist of dependent
variables, reflecting the fact that $\vv y \in \cM$ introduces
major difficulties in attempting to describe the measure theoretic
structure of either set.  This is true even in the specific case
that $\cM$ is a planar curve -- the main subject of this article.

\medskip

 In order to make any
reasonable progress  it is not unreasonable to assume that the
manifolds $ \cM$ under consideration are {\bf non-degenerate}.
Essentially, these are smooth sub-manifolds of $\R^n$ which are
sufficiently curved so as to deviate from any hyperplane.
Formally, a  manifold $\cM$ of dimension $m$ embedded in $\R^n$ is
said to be non-degenerate if it arises from a non--degenerate map
$\f:U\to \R^n$ where $U$ is an open subset of $\R^m$ and
$\cM:=\f(U)$. The map $\f:U\to \R^n:\bu\mapsto
\f(\bu)=(f_1(\bu),\dots,f_n(\bu))$ is said to be
\emph{non--degenerate at} $\bu\in U$ if there exists some $l\in\N$
such that $\f$  is $l$ times continuously differentiable on some
sufficiently small ball centred at $\bu$ and the partial
derivatives of $\f$ at $\bu$ of orders up to $l$ span $\R^n$. The
map $\f$ is \emph{non--degenerate} if it is non--degenerate at
almost every (in terms of $m$--dimensional Lebesgue measure) point
in $U$; in turn the manifold $\cM=\f(U)$ is also said to be
non--degenerate.  Any real, connected analytic manifold not
contained in any hyperplane of $\R^n$  is non--degenerate.

\vspace{7pt}

 Trivially, if the dimension $\dim \cM$ of the manifold
$\cM$ in $\R^n$ is strictly less than $n$ then $|\cM|_{\R^n} = 0$.
Thus, in attempting to develop a Lebesgue theory for the sets $ \cM
\cap \cS_n(\psi_1,\dots,\psi_n) $ and $ \cM \cap \cSM_n(\psi) $ it
is natural to use the induced Lebesgue measure $| \ . \ |_{\cM}$  on
$\cM$.

\vspace{5pt}

 In 1998,   D.\,Kleinbock \&
G.\,Margulis \cite{KM98} proved the Baker-Sprind\v{z}uk
conjecture:

\begin{theoremKM}\label{KM98}
Let $\cM$ be a non-degenerate manifold in $\R^n$. Then
\begin{equation}
\label{e:001} |\cM \cap \cSM_n(v)|_{\cM} \ = \ 0 \ \hspace{7mm}
{\rm if } \ \ \ \ v > 1 \ \  .
\end{equation}
\end{theoremKM}

By inclusion   (\ref{s:005}), Theorem~KM implies that for any
non-degenerate manifold
\begin{equation}
\label{sp} |\cM \cap \cS_n(v_1,\ldots,v_n) |_{\cM} \ = \ 0 \
\hspace{7mm} {\rm if } \ \ \ \ v_1+ \ldots +v_n > 1 \   . \ \ \ \
\footnote{When $v_1=\ldots=v_n$ this statement verifies the
conjecture of Sprind\v{z}uk which was  stated for analytic
manifolds only.} \end{equation}

\noindent Also, note that in view of (\ref{sv:004}) and
(\ref{sv:005}) both (\ref{e:001}) and  (\ref{sp}) are sharp.
 The first significant `clear cut' statement was for planar
curves. In 1964,  Schmidt \cite{Sch64}   established   (\ref{sp})
in the case that $\cM$ is a  $C^{(3)}$ non-degenerate planar curve
and $v_1= v_2$.

The result of Kleinbock \&{} Margulis gives some hope of developing
a general  metric theory for simultaneous approximation restricted
to manifolds,  analogous to that described in \S\ref{background}. As
stepping stones, it is natural to consider the following explicit
problems which ask for refinements of the measure zero statement of
Kleinbock \&{} Margulis.

\vspace{7pt}

\noindent{\bf Problem S1\,:} Given a non-degenerate manifold
$\cM\subset\R^n$ and $v>1 $ (respectively  $v_1+ \ldots +v_n >
1$), what is the Hausdorff dimension of $\cM \cap \cSM_n(v)$
(respectively $\cM \cap \cS_n(v_1,\ldots,v_n) \, $)  ?

\vspace{3pt}

\noindent{\bf Problem S2\,:} Given a non-degenerate manifold
$\cM\subset\R^n$ and an approximating function $\psi$
(respectively $\psi_1, \dots , \psi_n$), what is the weakest
condition  under which   $\cM \cap \cSM(\psi)$ (respectively $\cM
\cap \cS_n(v_1,\ldots,v_n)$ ) is of Lebesgue measure  zero?

\vspace{7pt}

 Problem S1 motivates the dimension theory for
simultaneous approximation restricted to manifolds whilst Problem
S2 motivates the convergent aspects of the Lebesgue theory. A
priori, convergent statements are usually easier to establish than
their divergent counterparts.

Until recently, the existing   metric theory  for simultaneous
approximation restricted to manifolds was  rather ad-hoc -- see
\cite{BDV02} for an account.  Even in the simplest geometric and
arithmetic  situation in which the manifold is a genuine curve in
$\R^2$  the above problems seemed to have been impenetrable.
However, in \cite{BDV02} we made significant progress towards
developing a complete metric theory for the sets $\cM \cap
\cS_2(\psi_1,\psi_2)$  with $\psi_1=\psi_2$ and $\cM$ a
non-degenerate planar curve\footnote{Note that in the case the
manifold  is a planar curve ${\cal C}$, a point on ${\cal C}$ is
non-degenerate if  the curvature at that point is non-zero. Thus,
${\cal C}$ is a non-degenerate planar curve if   the set of points
on ${\cal C}$ at which the curvature vanishes is a set of
one--dimensional Lebesgue measure zero. Moreover, it is not
difficult to show that the set  of points on a planar curve at
which the curvature vanishes but the curve  is non-degenerate is
at most countable. In view of this, the curvature completely
describes the non-degeneracy of planar curves.}. In this paper we
study the general simultaneous settings given to us  by the above
problems. This therefore includes the multiplicative setup. As in
\cite{BDV02}, we will mainly direct  our efforts towards the case
that the manifold $\cM$ is a  planar curve ${\cal C}$. Thus, $\dim
\cM=1$ and $n=2$ in the above problems.

\subsection{Statement of results}

\subsubsection{The  Lebesgue theory}

\begin{theorem}\label{t1}
Let $\psi_1,\psi_2$ be approximating  functions and let $\cC$ be a
$C^{(3)}$ non-degenerate planar curve. Then
$$
|\cC \cap \cS_2(\psi_1,\psi_2)|_{\cC} \, = \,    \mbox{\rm
F{\scriptsize ULL}} \qquad\text{{\rm
if}}\qquad\sum_{h=1}^\infty\psi_1(h)\psi_2(h)=\infty \ .
$$
\end{theorem}

The next theorem shows that the above result is best possible. We
establish the complementary `convergence result' for a class $\gQ$
of non-degenerate rational quadrics. A planar  curve $\cC$ is in
$\gQ$ if it is the image of either  the unit circle
$\cC_1:=\{(x_1,x_2) \in \R^2 :x_1^2+x_2^2=1\}$, the parabola
$\{(x_1,x_2) \in \R^2 : x_2 = x_1^2 \}$ or the hyperbola
$\{(x_1,x_2) \in \R^2 : x_1^2 - x_2^2 =1  \}$ under a rational
affine transformation of the plane.

\begin{theorem}\label{t2}
Let $\psi_1,\psi_2$ be approximating  functions and $\cC\in\gQ$.
Then
\begin{equation}\label{e:006}
\big|\cC \cap \cS_2(\psi_1,\psi_2)\big|_{\cC} \, = \,
\left\{\begin{array}{ll}
\mbox{\rm Z{\scriptsize ERO}}  & {\rm if} \;\;\;
\textstyle{\sum_{h=1}^\infty} \; \psi_1(h)\psi_2(h) \; <\infty \; ,\\
&
\\ \mbox{\rm F{\scriptsize ULL}} & {\rm if} \;\;\;
\textstyle{ \sum_{h=1}^\infty } \;  \psi_1(h)\psi_2(h) \;
 =\infty \; .
\end{array}\right.
\end{equation}
\end{theorem}

These theorems are a generalization of the results in \cite{BDV02}
which deal with the situation $\psi_1 = \psi_2$. The next theorem is
concerned with the multiplicative Lebesgue theory and is a
refinement of Theorem KM for manifolds in $\gQ$.

\begin{theorem}\label{t3}
Let $\psi$ be  an approximating  function  and $\cC\in\gQ$. Then
\begin{equation}\label{e:007}
|\cC \cap \cSM_2(\psi)|_{\cC} \, = \,  0 \qquad\text{{\rm
if}}\qquad\sum_{h=1}^\infty\psi(h)\log h<\infty \ .
\end{equation}
\end{theorem}

\subsubsection{The dimension theory }

Regarding Problem S1, for planar curves   we are able to give a
complete description for either form of simultaneous
approximation.

\begin{theorem}\label{t4}
Let $f\in C^{(3)}(I_0)$, where $I_0$ is an interval and $\cC_f :
=\set{(x,f(x)):x\in I_0}$. Let $v_1$ and $v_2$ be positive numbers
such that $0<\min(v_1,v_2)<1$ and $v_1+v_2\geq 1$.  Assume that
\begin{equation*} \dim\set{x\in I_0:
f''(x)=0}\le \frac{2-\min(v_1,v_2)}{1+\max(v_1,v_2)} \ .
\label{dimcond} \end{equation*} Then
\begin{equation}\label{e:008}
\dim\cC_f\cap\cS_2(v_1,v_2)\ = \
\frac{2-\min(v_1,v_2)}{1+\max(v_1,v_2)}\,.
\end{equation}
\end{theorem}


In the case $v_1=v_2$,  this theorem generalizes the dimension
results of \cite{BDV02}. Our next result is a general
$n$--dimensional statement concerning Lipshitz manifolds; i.e.
manifolds for which there exists an atlas of Lipshitz maps.

\begin{theorem}\label{t5}
Let $\cM$ be an arbitrary Lipshitz manifold in $\R^n$ of dimension
$\dim \cM$. Then
\begin{equation}\label{e:009}
\dim \cM\cap\cSM_n(v)\ \ge \  \dim \cM \, - \, 1 \ + \
\dfrac{2}{1+v} \hspace{9mm}  {\rm  if } \hspace{9mm} v \geq 1 \ \
.
\end{equation}
\end{theorem}

In essence, the above theorem  indicates  that the lower bound for
the Hausdorff dimension in the general multiplicative setup reduces
to a one dimensional problem. For clarification of this remark, see
(\ref{vb}) in \S\ref{pft5}.

\medskip

We conjecture that for manifolds in $\R^n$ which are  non-degenerate
everywhere except possibly on a set of dimension at most $\dim
\cM-1+ 2/(1+v)$, the lower bound given by Theorem~\ref{t5} is in
fact exact.  The following result verifies the conjecture for planar
curves.

\begin{theorem}\label{t6}
Let $f\in C^{(3)}(I_0)$, where $I_0$ is an interval and $\cC_f :
=\set{(x,f(x)):x\in I_0}$. Let $v > 1 $ and assume that
$\dim\set{x\in I_0: f''(x)=0}\le 2/(1+v)$. Then
\begin{equation}\label{e:010}
\dim\cC_f\cap\cSM_2(v) \ =  \ \dfrac{2}{1+v}  \ \ .
\end{equation}
\end{theorem}

\vspace{3ex}

\noindent {\em Remark.}
 Let $\psi$ be an approximating function
for which the limit
$$
\lambda (\psi) \ :=  \ \lim_{h \to \infty } \frac{ - \log \psi(h)
}{\log h }
$$
exists and is positive.  The quantity $\lambda (\psi)$ is usually
referred to as the {\em order} of $1/\psi$ and indicates the
limiting behavior of the function $1/\psi$ at  infinity.   On
making use of the fact that for any $\ve > 0 $, \begin{equation}
\label{joke}
 h^{- \lambda(\psi) - \ve } \ \leq \  \psi(h) \ \leq \
h^{- \lambda(\psi) + \ve } \end{equation} for all sufficiently
large $h$, the above dimension results (Theorems \ref{t4} --
\ref{t6}) can be easily generalized to approximating functions
$\psi$ for which the order $\lambda(\psi)$ exists.  For example,
Theorem \ref{t6} becomes

\begin{thsix}\label{t6*}
Let $f\in C^{(3)}(I_0)$, where $I_0$ is an interval and $\cC_f :
=\set{(x,f(x)):x\in I_0}$. Let $\psi$ be an approximating function
of order  $\lambda (\psi)
> 1 $ and assume that $\dim\set{x\in I_0: f''(x)=0}\le 2/(1+\lambda (\psi))$.
Then
\begin{equation*}\label{e:010T}
\dim\cC_f\cap\cSM_2(\psi) \ =  \ \dfrac{2}{1+\lambda (\psi)}  \ \
.
\end{equation*}
\end{thsix}

{\em Proof. }   Let $v := \lambda(\psi) $ and fix $\ve > 0 $ such
that $v-\ve > 1 $. In view of (\ref{joke}), it follows that
$$
\cC_f\cap\cSM_2(v+\ve)  \ \subset  \  \cC_f\cap\cSM_2(\psi) \
\subset \ \cC_f\cap\cSM_2(v - \ve)  \ .
$$
Theorem 6*  now follows from these inclusions, (\ref{e:009}) and
(\ref{e:010}), by letting $\ve \to 0 $.

\section{Preliminaries}

First some useful notation. For any point $\vv r\in\Q^n$ there
exists a smallest $q\in\N$ such that $q\vv r\in\Z^n$. Thus, every
point $\vv r\in\Q^n$ has a unique representation in the form
$$\frac{\bp}{q}=\frac{(p_1,\dots,p_n)}{q}=
\left(\frac{p_1}{q},\dots,\frac{p_n}{q}\right)$$ with
$(p_1,\dots,p_n)\in\Z^n$. Henceforth, we will only consider points
of $\Q^n$ in this form.  As usual, $C^{(n)}(I)$ will denote the set
of $n$--times continuously differentiable functions defined on some
interval $I$ of $\R$. Also, as usual   the Vinogradov symbols $\ll$
and $\gg$ will be used to indicate an inequality with an unspecified
positive multiplicative constant. If $a \ll b$ and $a \gg b$, we
write $a\asymp b$  and say that the quantities $a$ and $b$ are
comparable.

\subsection{Rational points close to a curve}

The following   estimates  on the number of rational points close
to a reasonably defined curve will be crucial towards establishing
our convergence  and (upper bound) dimension results.

\vspace{6pt}

 Let $I_0$ denote a finite, open interval of $\R$ and
let $f$ be a function in $C^{(3)}(I_0)$  such that
\begin{equation}\label{e:011}
\mbox{$c_1<|f''(x)|<c_2$ \ \ \ \ for all  \ $x\in I_0$} \ .
\end{equation}
%
%
%
%
Here $c_1$ and $  c_2$ are positive constants. Given an
approximating function $\psi$ and $Q \in \R^+$ consider the
counting function $N_f(Q,\p,I_0)$ given by $$ N_f(Q,\p,I_0) \ := \
\#\{\bp/q\in\Q^2 \, : \, q\le Q,\, p_1/q\in I_0,\,
|f(p_1/q)-p_2/q|<\p(Q)/Q\}. $$

\noindent In short, the function $N_f(Q,\p,I_0)$ counts the number
of rational points  with bounded  denominator lying within a
specified  neighbourhood of the curve $\cC_f := \{(x, f(x)): x \in
I_0 \}  $ parameterized by $f$. Now let
\begin{equation}\label{e:012}
\lim_{t\to+\infty}\psi(t)=\lim_{t\to+\infty}\ \frac{1}{t\psi(t)}=0.
\end{equation}

In \cite{Hux96}, Huxley obtains the following upper bound: {\em
For $\ve>0$ and $Q$ sufficiently large }
\begin{equation}\label{e:013}
N_f(Q,\p,I_0) \ \le \  \p(Q) \, Q^{2+\ve}.
\end{equation}
For this exact form of Huxley's estimate we refer the reader to
\cite[\S1.4]{BDV02}. In the case that the curve  is the unit
circle the above estimate can be sharpened. For $n \in \N $, let
$r(n)$ denote the number of representations of $n$ as the sum of
two squares.  A simple consequence of Theorem A in \cite[\S
A.1]{BDV02} is the following statement.

{\it There is a constant $C>0$ such that for any choice of real
numbers $Q$ and $\Psi$ satisfying
\begin{equation}\label{e:016S}
Q^{-1} \ (\log Q)^{260} \ \le \  \Psi \ < \ 1\qquad
\text{and}\qquad Q>1 \ \ ,
\end{equation}

one has that
\begin{equation}\label{e:017S}
 \sum_{Q < q \le
2Q } \ \ \ \sum_{\substack{n: \\ |q - \sqrt{n}| < \Psi }} \!\!\!\!
r(n) \ \, \le \, \ C\ \,  \Psi \ Q^2 \ \  .
\end{equation}
}

Notice that if $n$  is  the sum of two square, say $n= p_1^2 +
p_2^2$ then the inequality $|q-\sqrt n|<\Psi$ appearing in
(\ref{e:017S}) implies that the rational point $(p_1/q,p_2/q)$
lies within a constant times  $ \Psi/Q$ neighbourhood of  the unit
circle. Thus, we obtain the following sharpening of Huxley's
estimate.

{\em If  $\, \cC_f$ is the unit circle and
$ \psi(q)\ge q^{-1}(\log q)^{260}\, $ 
 for  all sufficiently large $q$, then
\begin{equation*}\label{e:014S}
    N_f(Q,\psi,I_0)\ll \psi(Q)Q^{2}  \ \ .
\end{equation*}
}

\noindent In fact, on adapting  the arguments of \cite[\S2]{BDV02}
it is relatively straightforward to extend the statement to any
planar curve $\, \cC_f$ in $\gQ$; i.e to any  non-degenerate
rational quadric. However, we shall not make use of this stronger
fact.

\subsection{Ubiquitous systems}

The divergence and (lower bound) dimension results stated in this
paper will be established via  a general technique developed in
\cite{BDV02}. The `general technique'  is based on the notion of
`ubiquity' as introduced in \cite{BDV03}.

\vspace{6pt}

Let $I_0$ be an interval in $\R$ and $\cR:=(\ra)_{\alpha\in\cJ}$ be
a family of {\em resonant points} $\ra$ of $I_0$ indexed by an
infinite set $\cJ$. Next let $\beta:\cJ\to \R^+:\alpha\mapsto\ba$ be
a positive function on $\cJ$.  Thus, the function $\beta$ attaches a
`weight' $\ba$  to the resonant point $\ra$. Also, for $t \in \N$
let $J_t:=\{\alpha \in \cJ:\ba\le 2^t\}$ and assume that  $\#J_t$ is
always finite.

Throughout, $\rho: \Rp \to\Rp$ will denote a function  satisfy
$\lim_{t\to\infty}\rho(t)=0$ and is usually referred to as the
{\em ubiquitous function}. Also $B(x,r)$ will denote the ball (or
rather the interval) centred at $x$ of radius $r$.

\begin{definition}[Ubiquitous systems on the real line]\sl
\label{US} Suppose there exists a ubiquitous function $\rho$ and
an absolute constant $\kappa > 0$ such that for any interval
$I\subseteq I_0$
\begin{equation*}
\liminf_{t\to\infty}\,\,\left|\,{\textstyle\bigcup_{\alpha\in
J_t}} \left(B(\ra,\rho(2^t)\right)\cap I)\right| \ \ge \  \kappa
\, |I| \ .
\end{equation*}
Then the system $(\cR;\beta)$ is called {\em locally ubiquitous in
$I_0$ with respect to $\rho$}.
\end{definition}

In \cite{BDV02} the theory of ubiquity is developed to incorporate
the situation in which the resonant points of interest lie within
some specified neighborhood of a given curve in $\R^n$.

With $n \geq 2$, let  $\cR:=(\ra)_{\alpha\in\cJ}$ be  a family of
{\em resonant points} $\ra$ of $\R^n$ indexed by an infinite set
$\cJ$. As before,  $\beta:\cJ\to \R^+:\alpha\mapsto\ba$ is  a
positive function on $\cJ$. For a point  $\ra$ in $\cR$,  let $\rak$
represent the $k$th coordinate of $\ra$.  Thus, $\ra :=
(R_{\alpha,1}, R_{\alpha,2}, \ldots, R_{\alpha,n})$. Throughout this
section and the remainder of the paper we will use the notation
$\cR_{\cC}(\Phi)$ to denote the sub-family of resonant points $\ra$
in $\cR$ which are ``$\Phi$--close'' to the curve $\cC=\cC_{\vv
f}:=\{(x,f_2(x),\dots,f_n(x)):x\in I_0\}$ where $\Phi$ is an
approximating function, $\vv f=(f_1,\dots,f_n):I_0\to\R^n$ is a
continuous map with $f_1(x)=x$ and $I_0$ is an interval in $\R$.
Formally, and more precisely  $$ \cRp := ( \ra )_{ \alpha\in
\cJ_\cC(\Phi) } \hspace{7mm} \text{where} \hspace{5mm}
\cJ_\cC(\Phi):= \{\alpha\in\cJ: \max\limits_{1\le k\le n}|
f_k(\rao)-\rak|<\Phi(\ba)\} \ . $$

\noindent Finally, we  will denote by $\cR_1$ the family of first
co-ordinates of the points in $\cRp $; that is $$ \cR_1 \ := \ (
R_{\alpha,1} )_{\alpha\in\cJp } \ \  . $$ \noindent By definition,
$\cR_1 $ is a subset of the interval $I_0$ and can therefore be
regarded as a set of resonant points for the theory of ubiquitous
systems in $\R$.  This leads us naturally to the following
definition in which the ubiquity function 
$\rho$ is as above.

\begin{definition}[Ubiquitous systems near curves]\sl \label{USNC}
The system $(\cRp,\beta)$ is called locally  ubiquitous with respect
to $\rho$ if the system $(\cR_1,\beta)$ is locally ubiquitous in
$I_0$ with respect to $\rho$.
\end{definition}

%
%

Next, given an approximating function $\Psi$ let
$\Lambda(\cRp,\beta,\Psi)$ denote the set  $x\in I_0$ for which
the system of inequalities
\[
\left\{\begin{array}{rcl} |x-\rao|&<&\Psi(\ba) \\[2ex]
\max\limits_{2\le k\le
n}|f_k(x)-R_{\alpha,k}|&<&\Psi(\ba)+\Phi(\ba) \ \
\end{array}\right.
\]
is satisfied  for infinitely many $\alpha\in\cJ$. The following
lemmas are stated and proved in \cite[\S3]{BDV02}.

\begin{lemma}\label{ktl}
Consider  the curve $\cC:=\set{(x,f_2(x),\dots,f_n(x)):x\in I_0}$,
where $f_2,\dots,f_n$ are locally Lipshitz in a finite interval
$I_0$.  Suppose that  $(\cRp,\beta)$ is a locally ubiquitous
system with respect to $\rho$. Let $\Psi$ be an approximating
function such that $\Psi(2^{t+1})\le \frac12 \Psi(2^t)$ for $t$
sufficiently large. Then
$$|\, \Lambda\left(\cRp,\beta,\Psi\right) \, | \ = \ |I_0| \ \ $$
whenever
$$
\sum_{t=1}^\infty\frac{\Psi(2^t)}{\rho(2^t)}=\infty
$$
\end{lemma}

\begin{lemma}\label{JTl}
Consider  the curve $\cC:=\set{(x,f_2(x),\dots,f_n(x)):x\in I_0}$,
where $f_2,\dots,f_n$ are locally Lipshitz in a finite interval
$I_0$.  Suppose that  $(\cRp,\beta)$ is a locally ubiquitous
system with respect to $\rho$ and let $\Psi$ be an approximating
function.
Then $$ \dim\Lambda \left(\cRp,\beta,\Psi \right) \ \ge \   d \,
:= \, \min\left\{1, \left| \limsup_{t\to\infty} \frac{\log
\rho(2^t)}{ \log \Psi(2^t)} \right| \right\}. $$
\end{lemma}

\begin{lemma}\label{cor1}
Let $I_0$ denote a finite, open interval of $\R$ and let $f$ be a
function in $C^{(3)}(I_0)$ satisfying\/ {\rm(\ref{e:011})}. Let
$\psi$ be an approximating function satisfying {\rm(\ref{e:012})}.
Let $\cC:=\{(x,f(x)) :x\in I_0\}$.   With reference to the
ubiquitous framework above, set
\begin{equation}\label{e:016}
\beta:\ \cJ := \Z^2 \times \N \to\N : (\bp,q) \to q  \ ,
\hspace{5mm} \Phi: t \to t^{-1}\psi(t) \hspace{5mm} \mbox{and}
\hspace{5mm} \rho :t \to u(t)/(t^2\psi(t))
\end{equation}
where $u:\Rp\to\Rp$ is any function such that
$\lim_{t\to\infty}u(t)=\infty$. Then the system
$(\Q^2_{\cC}(\Phi),\beta)$ is  locally  ubiquitous with respect to
$\rho$. 
\end{lemma}

\noindent{\em Remark.} In Lemma \ref{cor1}, the curve $\cC$ is
obviously a planar curve. Also, given $\alpha=(\bp,q) \in \cJ$  the
associated resonant point $\ra$ in the  ubiquitous system is simply
the rational point $ \bp/q$ in the plane.  Furthermore, $\cR :=
\Q^2$.

\section{Proof of Theorem~\ref{t1}}

As $\cC:=\cC_f$ is non-degenerate almost everywhere, we can restrict
our attention  to a sufficiently small patch of $\cC$, which can be
written as $\{(x,f(x)):x\in I\}$ where  $I$ is a sub-interval of
$I_0$  and $f $ satisfies  (\ref{e:011}) with $I_0$ replaced by $I$.
However, without loss of generality and for clarity, we assume that
$f $ satisfies (\ref{e:011}) on $I_0$.

%
%
%
%

We are given  that $\psi_1$ and $\psi_2$ are
 approximating functions such that
\begin{equation}\label{e:017}
\sum_{h=1}^\infty\psi_1(h)\psi_2(h)=\infty \ .
\end{equation}
Thus, at least one of the following two sums diverges:
$$
\sum_{h\in\N,\,\psi_1(h)\ge\psi_2(h)} \!\!\!\!\!\!\!\!\!\!\!\!
\psi_1(h)\psi_2(h) \qquad \qquad
\sum_{h\in\N,\,\psi_1(h)\le\psi_2(h)} \!\!\!\!\!\!\!\!\!\!\!\!
\psi_1(h)\psi_2(h) \ \ .
$$
Throughout, let us assume that the sum on the right is divergent.
The  argument below  can  easily be  modified to deal with the
case that only the sum on the left is divergent.


{\bf Step 1. \ } We show that there is no loss of generality in
assuming that
\begin{equation}\label{s:025}
\psi_2(h)\ge\psi_1(h)\qquad\text{for all}\quad h\in\N \ .
\end{equation}

\noindent  Define the auxiliary  function $\psi_1^*: h \to
\psi_1^*(h):=\min\{\psi_1(h),\psi_2(h)\}$. Then the sum
$$
\sum_{h=1}^\infty\psi_1^*(h)\psi_2(h)
$$
diverges since  by assumption it contains a divergent sub-sum. It
is readily verified that $\psi_1^*$ is an approximating function
and that $\cS_2(\psi_1^*,\psi_2)\subset\cS_2(\psi_1,\psi_2)$. Thus
to complete the proof of Theorem~\ref{t1} it suffices to prove the
result with $\psi_1$ replaced by $\psi_1^*$. Hence, without loss
of generality, (\ref{s:025}) can be assumed.


{\bf Step 2. \ } We show that there is no loss of generality in
assuming that
\begin{equation}\label{e:018}
  \psi_i(h)\to0\qquad\text{as}\qquad h\to\infty  \ \ \ \ \ \  \ (i=1,2) \ .
\end{equation}
 Define the increasing function $v : \R^+ \to \R^+$ as follows
$$
v(h):=\sum_{t=1}^{[h]}\psi_1(t)\psi_2(t) \ .
$$
In view of  (\ref{e:017}), $\lim_{t\to\infty}v(t)=\infty$. Fix
$k\in\N$. Then
$$
\sum_{t=k}^m\frac{\psi_1(t)\psi_2(t)}{v(t)} \ \ge \
\sum_{t=k}^m\frac{\psi_1(t)\psi_2(t)}{v(m)}\ = \
\frac{v(m)-v(k-1)}{v(m)} \ \to \ 1 \ \ \text{ as } \ \ m\to\infty.
$$
Hence
$$
\sum_{t=k}^\infty\frac{\psi_1(t)\psi_2(t)}{v(t)} \ \ge \ 1 \ \
\text{ for all }k.
$$
This implies that the sum $\sum_{t=1}^\infty
\psi_1(t)\psi_2(t)/v(t)$ diverges. Next, for $i=1,2$ consider  the
functions $$ \psi^*_i : h \to \psi^*_i(h):=\psi_i(h)/\sqrt{v(h)} \
.
$$  Then both $\psi^*_1(h)$ and $\psi^*_2(h)$ are
decreasing,  tend to $0$ as $h\to\infty$ and $\sum_{q=1}^\infty
\psi^*_1(q)\psi^*_2(q)=\infty$. Furthermore
$\cS_2(\psi^*_1,\psi^*_2)\subset\cS_2(\psi_1,\psi_2)$. Therefore,
it suffices to establish  Theorem~\ref{t1} for
$\psi^*_1,\psi^*_2$.


{\bf Step 3. \ } We show that there is no loss of generality in
assuming that
\begin{equation}\label{e:019}
\psi_2(h)\ge h^{-2/3} \ \text{ for all $h$.}
\end{equation}
To this end, define $\hat\psi_2(h)=\max\{\psi_2(h),h^{-2/3}\}$. In
view of (\ref{s:025}), it is readily verified that
\begin{eqnarray*}
\cS_2(\psi_1,\hat\psi_2) 
& \subseteq & \cS_2(\psi_1,\psi_2)\cup\cS_2(h\mapsto
h^{-2/3},h\mapsto h^{-2/3}) \ .
\end{eqnarray*}
By Schmidt's theorem  \cite{Sch64}, for almost all $x\in I_0$ we
have that $$ (x,f(x))\not\in\cS_2(h\mapsto h^{-2/3},h\mapsto
h^{-2/3})   \ .  $$ Hence
$$
\big|\{x\in I_0:(x,f(x))\in\cS_2(\hat\psi_1,\psi_2)\}\big| \ \leq
\  \big| \{x\in I_0:(x,f(x))\in\cS_2(\psi_1,\psi_2)\}\big| \ ,
$$
and to complete the proof of Theorem~\ref{t1}  it suffices to
prove that the set on the left has full measure. In turn, this
justifies (\ref{e:019}).


{\bf Step 4. \ } In view of Steps 2 and 3 above, the function
$\psi_2$ satisfies (\ref{e:012}) and  Lemma~\ref{cor1} is
applicable with  $\psi =\psi_2$. By (\ref{e:017}) and the fact
that $\psi_1$ and $\psi_2$ are decreasing we obtain that
$$
\infty=\sum_{t=0}^\infty\sum_{2^t\le h<2^{t+1}} \!\!\!\!\!\!
\psi_1(h)\psi_2(h) \  \le \ \sum_{t=0}^\infty\sum_{2^t\le
h<2^{t+1}} \!\!\!\!\!\! \psi_1(2^t)\psi_2(2^t) \ = \
\sum_{t=0}^\infty2^t \, \psi_1(2^t)\psi_2(2^t) \  .
$$
Hence
\begin{equation}\label{e:020}
\sum_{t=0}^\infty2^t \, \psi_1(2^t)\psi_2(2^t)=\infty \ .
\end{equation}

\noindent Next, define the increasing function $u:\R^+ \to \R^+ $
as follows
$$
u(h)=\sum_{t=0}^{[h]}2^t \, \psi_1(2^t)\psi_2(2^t) \ .
$$
Trivially,  $\lim_{t\to\infty}u(t)=\infty$. On using the same
argument as in Step 2 above we verify that
\begin{equation}\label{e:021}
\sum_{t=0}^\infty\frac{2^t \, \psi_1(2^t)\psi_2(2^t)}{u(t)} \ = \
\infty \ .
\end{equation}

Now let $\Phi(t):=\psi_2(t)/t$ and
$\rho(t):=u(\log_2t)/(t^2\psi_2(t))$.   By Lemma~\ref{cor1},
$(\Q^2_{\cC}(\Phi),\beta)$ is locally ubiquitous relative to
$\rho$, where $\beta$ is given by (\ref{e:016}). Let
$\Psi(t)=\psi_1(t)/t$. In view of (\ref{e:021}),
$$
\sum_{t=1}^\infty\frac{\Psi(2^t)}{\rho(2^t)} \ = \
\sum_{t=1}^\infty\frac{\quad\displaystyle
\frac{\psi_2(2^t)}{2^t}\quad}{\displaystyle\frac{u(t)}{2^{2t}\psi_1(2^t)}}
\  = \ \sum_{t=1}^\infty\frac{2^t\psi_1(2^t)\psi_2(2^t)}{u(t)} \ =
\ \infty \ .
$$
Since $\psi_1$ is decreasing,
$$
\Psi(2^{t+1}) \ := \ \frac{\psi_1(2^{t+1})}{2^{t+1}} \ \le \
\frac12\cdot\frac{\psi_1(2^t)}{2^t} \ := \ \frac12  \Psi(2^t) \ .
$$
Thus the conditions of  Lemma~\ref{ktl} are satisfied and it
follows that $\Lambda(\Q^2_{\cC}(\Phi),\beta,\Psi)$ is of  full
measure. By definition  and (\ref{s:025}), the  set
$\Lambda(\Q^2_{\cC}(\Phi),\beta,\Psi)$ consists of points $x \in
I_0$ such that the system
$$
\left\{
\begin{array}{l}
  \big|x-\frac{p_1}{q}\big|<\Psi(q)=
\frac{\psi_1(q)}{q}<\frac{2\psi_1(q)}{q},  \\[1ex]
\big|f(x)-\frac{p_2}{q}\big|<\Psi(q)+\Phi(q)=\frac{\psi_1(q)}{q}+
\frac{\psi_2(q)}{q}\le
\frac{2\psi_2(q)}{q}
\end{array}
\right.
$$
has infinitely many solutions  $\bp/q\in\Q^2$. Obviously for
$x\in\Lambda(\Q^2_{\cC}(\Phi),\beta,\Psi)$ the point $(x,f(x))$ is
in $\cS_2(2\psi_1,2\psi_2)$. To complete the proof of
Theorem~\ref{t1} we simply apply what has  already been proved to
the approximating functions $\frac12\psi_1$ and $\frac12\psi_2$.
\vspace{-0ex}

\hfill $ \spadesuit$

\section{Proof of Theorem \ref{t3}}

By definition any rational quadric $\cC\in\cQ$ is the image of
either the unit circle $\cC_1:=\{(x_1,x_2) \in \R^2
:x_1^2+x_2^2=1\}$, the parabola $\{(x_1,x_2) \in \R^2 : x_2 =
x_1^2 \}$ or the hyperbola $\{(x_1,x_2) \in \R^2 : x_1^2 - x_2^2
=1  \}$ under a rational affine transformation of the plane. It is
easily verified that the measure `zero' statement of  Theorem
\ref{t3} is invariant under rational affine transformations of the
plane. In view of this, it suffices to establish the statement of
the theorem for  the unit circle, the parabola and the hyperbola.
Below, we only consider the case of the unit circle $\cC_1$ and
leave the hyperbola and parabola to the reader. The required
modifications are relatively straightforward once the reader is
armed with the arguments appearing in \cite[\S2.1]{BDV02}.

\vspace{2ex}

From this point onwards $\cC = \cC_1$ -- the unit circle. First,
notice that it suffices to prove the theorem for every arc of
$\cC_1$ given by $\cC_1^\ve=\{(x,y)\in\cC_1:\ve<x,y<1-\ve\}$ with
$\ve>0$. Next, we are given that
\begin{equation} \label{dmsum}
\sum_{q=1}^\infty\psi(q)\log q<\infty  \ .
\end{equation}

\noindent Therefore, without loss of generality we can assume that
\begin{equation}\label{eS:025}
q^{-1} (\log q)^{-3}  \ <  \ \psi(q) \ < \ q^{-1} (\log q)^{-1}
\quad\text{for  sufficiently large } q  \ .
\end{equation}

\vspace{2ex}

\noindent To see this, note that since $\psi$ is decreasing
$$
\sum_{h/2<q<h}\psi(q)\log q \ \ge \  \sum_{h/2<q<h}\psi(h)\log
(h/2) \ \asymp \  h \, \psi(h) \; \log h \
$$
for every natural number $h$. In view of (\ref{dmsum}), we have
that  the left hand side of the above inequality tends to zero as
$ h \to \infty$. It follows that
$$
h \, \psi(h)  \, \log h  \to 0 \quad\text{as } \ h \to \infty \ \
.
$$
This establishes the right hand side inequality of (\ref{eS:025}).
Next, if  the left hand side inequality of (\ref{eS:025}) is not
satisfied then we replace $\psi$ with the auxiliary  function
$$\widetilde{\psi} : q \to \widetilde{\psi}(q) :=
\max\{\psi(q),\ q^{-1} (\log q)^{-3}
 \}  \ . $$ Then $\widetilde{\psi}$ is clearly an approximating  function
for which  both (\ref{dmsum}) and   the left hand side inequality
of (\ref{eS:025}) are satisfied  with $\psi $ replaced by
$\widetilde{\psi} $. Furthermore, $$\cSM_2(\widetilde\psi) \
\supset \ \cSM_2(\psi) \ .  $$ Thus it suffices to prove the
theorem with $\psi$ replaced by $\widetilde\psi$. Hence, without
loss of generality, (\ref{eS:025}) can be assumed.

\vspace{2ex}

 The $\limsup$ set $\cC_1^\ve\cap\cSM_2(\psi)$ has the
following natural representation:
$$
  \cC_1^\ve\cap\cSM_2(\psi)=\bigcap_{n=1}^\infty \ \bigcup_{q=n}^\infty\
\bigcup_{(p_1,p_2)\in\;\Z^2}\
\Big\{(x,y)\in\cC_1^\ve:\Big|x-\frac{p_1}{q}\Big|
\cdot\Big|y-\frac{p_2}{q}\Big|<\frac{\psi(q)}{q^2}\Big\}.
$$
Using the fact that  $\psi$ is decreasing, we have that for any
$n$
\begin{equation}\label{eS:026}
\cC_1^\ve\cap\cSM_2(\psi)\subset \bigcup_{t=n}^\infty\
\bigcup_{2^t\le q< 2^{t+1}}\ \bigcup_{(p_1,p_2)\in\;\Z^2}\
\Big\{(x,y)\in\cC_1^\ve:\Big|x-\frac{p_1}{q}
\Big|\cdot\Big|y-\frac{p_2}{q}\Big|<\frac{\psi(2^t)}{(2^t)^2}\Big\}.
\end{equation}

\noindent If $t\in\N$,  $(x,y)\in\cC_1^\ve$, $q\in\N$ with $2^t\le
q<2^{t+1}$ and
$\Big|x-\frac{p_1}{q}\Big|\cdot\Big|y-\frac{p_2}{q}\Big|<\frac{\psi(2^t)}{(2^t)^2}$
for some $(p_1,p_2)\in\Z^2$, then  there is a unique integer $m$
such that
$$
2^{m-1}\,\frac{\sqrt{2\,\psi(2^t)}}{2^t}\le
\Big|x-\frac{p_1}{q}\Big|<2^m\,\frac{\sqrt{2\,\psi(2^t)}}{2^t}.
$$
Therefore for this number $m$ we also have that
$$
\Big|y-\frac{p_2}{q}\Big|< 2^{-m}\frac{\sqrt{ 2\,\psi(2^t)}}{2^t}
\  .
$$
The upshot of this is that
\begin{equation}\label{eS:027}
\cC_1^\ve\cap\cSM_2(\psi)\subset \bigcup_{t=n}^\infty\
\bigcup_{2^t\le q< 2^{t+1}}\ \bigcup_{(p_1,p_2)\in\;\Z^2}\
\bigcup_{m=-\infty}^{+\infty} \cC_1^\ve\cap S(q,p_1,p_2,m) \ ,
\end{equation}
where
$$
S(q,p_1,p_2,m)=\Big\{(x,y)\in\R^2:\Big|x-\frac{p_1}{q}\Big|<
2^{m}\frac{\sqrt{2\,\psi(2^t)}}{2^t},\ \ \
\Big|y-\frac{p_2}{q}\Big|<2^{-m}\frac{\sqrt{
2\,\psi(2^t)}}{2^t}\Big\}
$$
and $t$ is uniquely defined by $2^t\le q< 2^{t+1}$. The aim is to
show that the Lebesgue measure $ | \ . \  |_{\cC_1^\ve} $ of the
R.H.S. of (\ref{eS:027}) tends to zero as $n \to \infty$. Since for
each $n$ the R.H.S. of (\ref{eS:027}) is a cover for $\cC_1^\ve \cap
\cSM_2(\psi)$, it follows that $|\cC_1^\ve \cap
\cSM_2(\psi)|_{\cC_1^\ve} =   0 $ as required.   To proceed, we
consider two cases. Namely,  case (a): $m \in \Z$ such that
\begin{equation} \label{casea}
2^{-|m|} \ \ge \  t\sqrt{\psi(2^t)} \  \ ~
\end{equation}
and case (b): $m \in \Z$ such that
\begin{equation} \label{caseb}
2^{-|m|} \ \le \  t\sqrt{\psi(2^t)} \  \ .
\end{equation}


{\bf Case (a) }: \ First,  observe that (\ref{casea}) together with
(\ref{eS:025}) implies that
\begin{equation}\label{eSV:028}
t \ \geq \  2 \,  |m| \ \ \ .
\end{equation}
Next, it is a simple mater to see that
\begin{equation}\label{eS:028}
|\,\cC_1^\ve\cap S(q,p_1,p_2,m)|_{\cC_1^\ve} \ \ll \  2^{-|m|}\
\frac{\sqrt{\psi(2^t)}}{2^t} \ 
\ .
\end{equation}
The implied constant depends on only  $\ve$ and  is therefore
irrelevant to the rest of the argument.

Given $t$ and $m$, let $N(t,m)$ denote the number of triples
$(q,p_1,p_2)$ with $2^t\le q< 2^{t+1}$ such that $\cC_1^\ve\cap
S(q,p_1,p_2,m)\not=\emptyset$.  Suppose that $\cC_1^\ve \cap
S(q,p_1,p_2,m) \neq \emptyset$.  Then for some $(x,y)\in\cC_1^\ve$
and  $\theta_1, \theta_2$ satisfying $-1<\theta_1,\theta_2<1$, we
have that
$$
x \ = \ \frac{p_1}{q} \, + \, \theta_1 \; 2^{|m|} \,
\frac{\sqrt{2\,\psi(2^t)}}{2^t} \ , \ \  \ \  \ \ y \ = \
\frac{p_1}{q} \, + \, \theta_2 \; 2^{|m|} \;
\frac{\sqrt{2\,\psi(2^t)}}{2^t} \ \ .
$$
Hence
\begin{eqnarray*}
1 & = &   x^2+y^2  \\[1ex]& = & \sum_{i=1}^2\left(\frac{p_i}{q} \, + \,
\theta_i \, 2^{|m|}\frac{\sqrt{2\,\psi(2^t)}}{2^t} \; \right)^2
 \\ & & \nonumber \\
 & = & \frac{p_1^2+p_2^2}{q^2} \ + \
\frac{p_1\theta_1+p_2\theta_2}{q} \ \  2^{|m|+1} \ \
\frac{\sqrt{2\,\psi(2^t)}}{2^t}
 \ + \ (\theta_1^2+\theta_2^2) \ 2^{2|m|+1} \ \frac{\psi(2^t)}{2^{2t}} \ .
\end{eqnarray*}
It follows that
$$
| \ q^2-p_1^2+p_2^2 \ | \ \ll \ q \ \max\{|p_1|,|p_2|\} \
2^{|m|-t} \ \sqrt{\psi(2^t)} \  \ + \  \ q^2 \, 2^{2|m|-2t} \,
\psi(2^t) \ \  .
$$
On dividing both sides of the inequality by $q+\sqrt{p_1^2+p_2^2}$
and using the fact that $\sqrt{p_1^2+p_2^2}\ge
\max\{|p_1|,|p_2|\}$ we obtain that
$$
| \ q- \textstyle{\sqrt{p_1^2+p_2^2}} \ | \ \ll \  2^t \, U_{m,t}
\; (1+U_{m,t}) \ ,
$$
where
$$
U_{m,t} \ := \  2^{-t} \, 2^{|m|} \, \sqrt{\psi(2^t)}  \
\stackrel{(\ref{casea})}{\ \leq \ } \ 2^{-t} \, t^{-1} \ < \ 1 \ .
$$
Thus
\begin{equation}\label{eS:029}
| \ q- \textstyle{ \sqrt{p_1^2+p_2^2}} \ | \ \ll \ 2^{|m|} \,
\sqrt{\psi(2^t)} \ \stackrel{(\ref{casea})}{\ \leq \ } \  t^{-1} \ .
\end{equation}
The upshot of this  is that there exists an absolute constant $c>0$
such that
$$
N(t,m)  \  \ll  \ \sum_{2^t\le q<2^{t+1}}\ \ \sum_{|q-\sqrt n|<
c2^{|m|}\sqrt{\psi(2^t)}} r(n)\ .
$$

Now set $Q=2^t$ and $\Psi=c2^{|m|}\sqrt{\psi(2^t)}$. In view of
(\ref{eS:025}) and (\ref{eS:029}) we have that (\ref{e:016S}) is
satisfied  for all sufficiently large $t$, independently of $m$.
Hence, (\ref{e:017S}) implies that
\begin{equation}\label{eS:030}
N(t,m) \ \ll \  2^{|m|} \, 2^{2t} \, \sqrt{\psi(2^t)}   \ ,
\end{equation}
where the implied constant is independent of both $t$ and $m$.

\medskip

It now follows,  via (\ref{eS:028})  and (\ref{eS:030}) that the
Lebesgue measure $ | \ . \ |_{\cC_1^\ve} $ of the R.H.S. of
(\ref{eS:027}) restricted to case (a) is bounded above by
\begin{eqnarray*}
\sum_{t=n}^\infty\ \sum_{m\in\text{Case \!\!(a)}} \!\!\! 2^{-|m|}\
\frac{\sqrt{\psi(2^t)}}{2^t} \times 2^{|m|} \,
2^{2t}\sqrt{\psi(2^t)} & \ll &  \sum_{t=n}^\infty\
\sum_{m\in\text{Case \!\!(a)}} \!\! 2^t \, \psi(2^t)
\\ & \stackrel{(\ref{eSV:028})}{\ \ll \ } &
\sum_{t=n}^\infty\ t\, 2^t\psi(2^t) \\ & \asymp &
\sum_{q=2^n}^\infty\
 \psi(q)\, \log q  \ \ .
\end{eqnarray*}
The above comparability follows from the fact that $\psi$ is an
approximating function and therefore decreasing.  In view of
(\ref{dmsum})
$$
\sum_{q=2^n}^\infty\
 \psi(q)\, \log q \ \ \to \ \ 0\qquad\text{as }\qquad n\to\infty  \ ,
$$
and so the Lebesgue measure $ | \ . \  |_{\cC_1^\ve} $ of the R.H.S.
of (\ref{eS:027}) restricted to case (a)  tends to zero as $n \to
\infty$.

{\bf Case (b) }:   \   In view of (\ref{caseb}), we have that
\begin{equation}
S(q,p_1,p_2,m) \ \subset \  S'(q,p_1)\times[0,1] \; \cup \;
[0,1]\times S'(q,p_2)
\end{equation}
where
\begin{equation}
S'(q,p)  \ = \ \Big\{s\in[0,1]:
\Big|s-\frac{p}{q}\Big|<\frac{2\,t\,\psi(2^t)}{2^t}\Big\}
\end{equation}
Thus, the R.H.S. of (\ref{eS:027}) restricted to case (b)  is
contained in the following set:
\begin{equation}\label{eS:033}
\bigcup_{t=n}^\infty\ \bigcup_{2^t\le q< 2^{t+1}}\
\bigcup_{(p_1,p_2)\in\;\Z^2}\ \cC_1^\ve \ \cap \ (\,
S'(q,p_1)\times[0,1]\cup[0,1]\times S'(q,p_2) \,  ) \  \ .
\end{equation}
It is readily verified that for any choice of $p_1,p_2$ and $ q $
appearing in (\ref{eS:033}),
$$
| \, \cC_1^\ve \ \cap \ (\, S'(q,p_1)\times[0,1]\cup[0,1]\times
S'(q,p_2) \, ) \, |_{\cC_{1}^{\ve}}  \ \ll \ t\,\psi(2^t)/2^t   \  \
.
$$
The implied constant depends only on $\ve$ and is therefore
irrelevant. Furthermore, for a fixed $t$ and $q$ in (\ref{eS:033})
the number of $(p_1,p_2)\in\;\Z^2 $  for which the sets  $$
\cC_1^\ve \ \cap \ (\, S'(q,p_1)\times[0,1]\cup[0,1]\times S'(q,p_2)
\, )
$$ are non-empty and disjoint  is $\ll q $.   It now follows that the
Lebesgue measure $ | \ . \  |_{\cC_1^\ve} $
 of the set given by (\ref{eS:033}) is bounded above by
 $$ \sum_{t=n}^\infty t\,\psi(2^t) \ \asymp  \
\sum_{q=2^n}^\infty \psi(q)\log q   \ \ \to \ \ 0\qquad\text{as }
n\to\infty  \  .
$$
Hence the Lebesgue measure $ | \ . \  |_{\cC_1^\ve} $ of the R.H.S.
of (\ref{eS:027}) restricted to case (b)  tends to zero as $n \to
\infty$.

\medskip

The upshot of cases (a) and (b) is that the Lebesgue measure $ | \ .
\  |_{\cC_1^\ve} $ of the R.H.S. of (\ref{eS:027})  tends to zero as
$n \to \infty$ and so
$$
|\cC_1^\ve \cap \cSM_2(\psi)|_{\cC_1^\ve}  \ =  \   0  \ \ .
$$
This completes the proof of Theorem \ref{t3}.

 \vspace{-0ex}

\hfill $ \spadesuit$

\section{Proof of Theorem~\ref{t2}}

The divergence part of Theorem \ref{t2} is a consequence of
Theorem~\ref{t1}. Thus we proceed with establishing the convergence
part of Theorem \ref{t2}. To a certain degree the proof of this
follows the same line of argument as the proof of  Theorem \ref{t3}.
In particular, it suffices to establish the convergent statement of
the theorem for the unit circle, the parabola and the hyperbola. As
in the proof of Theorem \ref{t3}, we consider the case of the unit
circle $\cC_1$ only  and leave the hyperbola and parabola to the
reader.

\medskip

Let $\cC = \cC_1$ -- the unit circle,  and  notice that it suffices
to prove the theorem for every arc of $\cC_1$ given by
$\cC_1^\ve=\{(x,y)\in\cC_1:\ve<x,y<1-\ve\}$ with $\ve>0$. For the
sake of convenience, let $\psi := \psi_1 $ and $\phi := \psi_2$. It
is clear that
$$
\cS_2(\psi,\phi) \ \subset \  \cS_2(\psi^*,\psi_*) \cup
\cS_2(\psi_*,\psi^*)
$$
where
$$
\psi_*=\min\{\psi,\phi\} \ ,  \ \ \ \ \ \ \psi^*=\max\{\psi,\phi\}.
$$
Since $\psi^*\psi_*=\psi\phi$, we have that $\sum\psi^*(q)\psi_*(q)
<  \infty$. Thus to prove the theorem it suffices to prove that both
the sets  $\cC_1^\ve \cap \cS_2(\psi^*,\psi_*)$ and $\cC_1^\ve \cap
\cS_2(\psi_*,\psi^*)$ are of  Lebesgue measure $ | \ . \
|_{\cC_1^\ve} $ zero. We will consider one of these two sets -- the
other case is similar. Thus, without loss of generality we assume
that
$$
\psi(q) \ \ge \ \phi(q)\text{ \ \ for all }q \ \ .
$$

\medskip

Since $\sum_{q=1}^\infty\psi(q)\, \phi(q)<\infty$ and both $\psi $
and $\phi$ are decreasing we have that $ \psi(q)\phi(q)<q^{-1} $ for
all sufficiently large $q$. Hence
$$
\phi(q) \ \le \  q^{-1/2}  \text{\qquad for sufficiently large } q \
.
$$
Further, we can assume that
\begin{equation} \label{th2big}
\psi(q)\ge q^{-2/3}\text{\qquad for all }q\in\N  \ \ .
\end{equation}
To see  this, consider the auxiliary  function $\widetilde{\psi} : q
\to \widetilde{\psi}(q) := \max\{\psi(q),q^{-2/3}\}$. Clearly,
$\widetilde{\psi}$ is an approximating function, satisfies
(\ref{th2big})  and
$$
\cS_2(\psi,\phi)\subset\cS_2(\widetilde{\psi},\phi) \ .
$$
Moreover,
$$
\begin{array}{rl}
\displaystyle\sum_{q=1}^\infty\widetilde{\psi}(q) \, \phi(q)\le&
\displaystyle \sum_{q=1}^\infty\psi(q)\phi(q) + \sum_{q=1}^\infty
q^{-2/3}\phi(q)
 \\[4ex]
\ll&\displaystyle\sum_{q=1}^\infty\psi(q)\phi(q) +
\sum_{q=1}^\infty q^{-2/3}q^{-1/2} \ < \ \infty \ \ .
\end{array}
$$
Thus, it suffices to prove the convergence part of Theorem \ref{t2}
with $\psi:=\psi_1$ replaced $\widetilde{\psi}$. Hence, without loss
of generality, (\ref{th2big}) can be assumed.

\vspace{2ex}

In analogy  to  (\ref{eS:026}), it is readily verified that for any
$n \geq 1 $
\begin{equation} \label{incl}
\cC^\ve_1\cap\cS_2(\psi,\phi) \ \subset \  \bigcup_{t=n}^\infty\
\bigcup_{2^t\le q< 2^{t+1}}\ \bigcup_{(p_1,p_2)\in\;\Z^2}\ \!\!
\cC^\ve_1\cap S_2(p_1,p_2,q) \ \ ,
\end{equation}
where
$$
S_2(p_1,p_2,q)=\Big\{(x,y)\in\R^2:\Big|x-\frac{p_1}{q}\Big|<\frac{\psi(2^t)}{2^t},\
\ \Big|y-\frac{p_2}{q}\Big|<\frac{\phi(2^t)}{2^t}\Big\}
$$
and $t$ is uniquely defined by $2^t\le q< 2^{t+1}$. Next,  in
analogy to  (\ref{eS:028}), we verify that
\begin{equation} \label{minq}
|\,\cC^\ve_1\cap S_2(q,p_1,p_2)|_{\cC^\ve_1}  \ \ll \
\frac{\phi(2^t)}{2^t} \ .
\end{equation}
Again, the implied constant depends only on $\ve$ and is therefore
irrelevant to the rest of the argument. For $t$ fixed , let $N(t)$
denote the number of triples $(q,p_1,p_2)$ with $2^t\le q< 2^{t+1}$
such that $\cC^\ve_1 \cap S(q,p_1,p_2)\not=\emptyset$. On modifying
the argument used to establish (\ref{eS:029}) and (\ref{eS:030}) in
the proof of Theorem \ref{t3},  one obtains that
\begin{equation} \label{cinq}
N (t) \ \ll  \  2^{2t}\psi(2^t) \ .
\end{equation}
It is worth stressing that the argument within the proof of Theorem
\ref{t3} is much  simplified in the  current situation due to the
absence of the additional parameter $m$.

\medskip

The upshot of the above inclusions and estimates is that
\begin{eqnarray*}
|\;\cC^\ve_1\cap\cS_2(\psi,\phi)|_{\cC^\ve_1}
&\stackrel{(\ref{incl})}{\ \ll \ } & \sum_{t=n}^\infty\ \sum_{2^t\le
q< 2^{t+1}}\ \bigcup_{(p_1,p_2)\in\;\Z^2}\ |\cC^\ve_1\cap
 S_2(q,p_1,p_2)|\\
 &\stackrel{(\ref{minq})}{\ \ll \ } & \sum_{t=n}^\infty\ N(t)\;
 \frac{\phi(2^t)}{2^t} \\ &  \stackrel{(\ref{cinq})}{\ \ll \ }&  \sum_{t=n}^\infty \ 2^t
 \psi(2^t)\phi(2^t)\\
 &\asymp& \sum_{q= 2^n}^\infty \  \psi(q)\phi(q)
\\
\end{eqnarray*}
Since $\sum_{q= 2^n}^\infty   \psi(q)\phi(q) < \infty $, we have
that $\sum_{q= 2^n}^\infty   \psi(q)\phi(q)\to 0 $ as $ n\to\infty
$.  Thus,
$$
|\;\cC^\ve_1\cap\cS_2(\psi,\phi)|_{\cC^\ve_1} = 0  \   .
$$
This completes the proof of the theorem. \vspace{-2ex}

\hfill $ \spadesuit$

\section{Proofs of Theorems \ref{t4} -- \ref{t6}}

%
%
%

\subsection{Proof of Theorem \ref{t4}}

The statement of the theorem will follow on establishing the upper
and lower bounds for the dimension separately. Without loss of
generality we can  assume that $f$ satisfies (\ref{e:011}) on
$I_0$ (see \cite[\S5]{BDV02} if necessary)  and that $v_1\ge v_2$.
In view of the latter, our aim is to show that
$$
\dim\cC_f\cap\cS_2(v_1,v_2) \ = \ \frac{2-v_2}{1+v_1} \ .
$$

\vspace{3ex}

\noindent{\em The upper  bound. \  } For a point $\bp/q\in\Q^2$,
define
$$
\sigma(\bp/q):=\{(x,y)\in \R^2:|x-p_1/q|<q^{-v_1-1},\
|y-p_2/q|<q^{-v_2-1}\}.
$$
In view of (\ref{e:011}) and the fact that $I_0$ is a bounded
interval we have that $f'$ is bounded on $I_0$ and so
$|\cC_f\cap\sigma(\bp/q)|\ll q^{-v_1-1}$. Clearly, if
$\sigma(\bp/q)\cap\cC_f\not=\emptyset$ then the distance of $\bp/q$
from $\cC_f$ is at most  a constant times $ q^{-1-v_2}$. Let
$\ve>0$. In view of (\ref{e:013}) we have  that for $t$ sufficiently
large the number $\bp/q \in \Q^2$ with $2^t\le q< 2^{t+1}$ and
$\sigma(\bp/q)\cap\cC_f\not=\emptyset$ is  at most
$2^{t(2+\ve-v_2)}$. Now let
$$\eta \ := \
\frac{2-v_2+2\ve}{1+v_1}\ .
$$
Then
\begin{eqnarray*}
\sum_{\bp/q\in\Q^2 \, : \,\cC_f\cap\sigma(\bp/q)\not=\emptyset}
\!\!\!\!\!\!\!\!\!\!\!\! \diam(\cC_f\cap\sigma(\bp/q))^\eta
 & = & \sum_{t=0}^\infty \ \ \ \ \sum_{\bp/q\in\Q^2,\,\cC_f\cap\sigma(\bp/q)
\not=\emptyset,\,2^t\le q<2^{t+1}}
\!\!\!\!\!\!\!\!\!\!\!\!\!\!\!\!\!\!\!\!\!\!
\diam(\cC_f\cap\sigma(\bp/q))^\eta
\\ \\ & \ll &
\sum_{t=0}^\infty 2^{t(2+\ve-v_2)}\cdot2^{t(-1-v_1)\eta} \ = \
\sum_{t=0}^\infty 2^{- t\ve} \ < \ \infty \ \ .
\end{eqnarray*}
By the Hausdorff--Cantelli Lemma \cite[p.~68]{BD99},
$\dim\cC_f\cap\cS_2(v_1,v_2)\le\eta $. As $\ve>0$ is arbitrary,
\begin{equation}\label{e:031}
\dim\cC_f\cap\cS_2(v_1,v_2) \ \le\ \frac{2-v_2}{1+v_1} \ .
\end{equation}

\vspace{3ex}

\noindent{\em The lower bound. \  }Firstly, with reference to
Lemma \ref{cor1} let $\psi(t):=\frac12t^{-v_2}$ and $u(t):=
t^{\ve} $ where $\ve >0 $ is arbitrary. Thus
$\Phi(t)=\frac12t^{-1-v_2}$ and $\rho(t):=\frac12t^{-2+v_2+\ve}$.
Since $0 < v_2 = \min(v_1,v_2) < 1 $, the approximating function
$\psi$ satisfies (\ref{e:012}) and it follows that
$(\Q^2_{\cC}(\Phi),\beta)$ is locally ubiquitous with respect to
$\rho$. Next, let  $\Phi(t):=\frac12t^{-1-v_1}$. Then
Lemma~\ref{JTl} implies that
$$
\dim\Lambda \left(\cRp,\beta,\Psi \right) \ \ge \ \min\left\{1,
\left| \limsup_{t\to\infty} \frac{\log
\frac122^{t(-(2-v_2-\ve))}}{ \log \frac122^{t(-(1+v_1))}} \right|
\right\}\ = \ \frac{2-v_2-\ve}{1+v_1}.
$$
As $\ve>0$ can be made arbitrarily small, we have that $\dim\Lambda
\left(\cRp,\beta,\Psi \right) \ge \frac{2-v_2}{1+v_1}$. Finally,  it
is readily verified  that
$$
\widetilde\Lambda
\left(\cRp,\beta,\Psi\right):=\{(x,f(x)):x\in\Lambda
\left(\cRp,\beta,\Psi\right)\} \  \subset \
\cC_f\cap\cS_2(v_1,v_2).
$$
Hence
$$
\dim \cC_f\cap\cS_2(v_1,v_2) \ \geq \  \dim\widetilde\Lambda
\left(\cRp,\beta,\Psi \right) \ = \  \dim\Lambda
\left(\cRp,\beta,\Psi \right)\  \geq \ \frac{2-v_2}{1+v_1} \ .
$$
The equality here is justified by the fact that the map
$x\mapsto(x,f(x))$ is locally bi-Lipshitz.


\hfill $ \spadesuit$

\subsection{Proof of Theorem~\ref{t5} \label{pft5}}

Let $m=\dim \cM$. Since $\cM$ is a  Lipshitz manifold in $\R^n$,
there exists  a local parameterization of $M$ of the form $\vv
f=(x_1,\ldots,x_m,f_{m+1},\dots,f_n)$  where $\vv f$ is an
invertible continuous map of $x_1,\ldots,x_m$ defined on $\R^m$
such that $\vv f^{-1}$ satisfies the Lipshitz condition. It is
easy to verify  that any point on $\cM$ with $x_1\in\cSM_1(v)$
belongs to $\cSM_n(v)$. Therefore,
$$
B:= \vv f(\cSM_1(v)\times\R^{m-1}) \ \subset \  \cSM_n(v)\cap \cM
\ \ .
$$
Since $\vv f^{-1}$ is a Lipshitz map and $\vv
f^{-1}(B)=\cSM_1(v)\times\R^{m-1}$, it follows  that for $v \geq
1$
\begin{eqnarray}
\dim \cSM_n(v)\cap \cM \ \ge \  \dim B \ \ge \ \dim
\cSM_1(v)\times\R^{m-1}& = & m-1+\dim\cSM_1(v) \label{vb} \\[2ex] & = &
\dim \cM-1+\frac{2}{1+v} \nonumber \ .
\end{eqnarray}

The fact that $\dim\cSM_1(v)= 2/(1+v)$  is   the
Jarn\'{\i}k--Besicovitch theorem (see \S\ref{background}).
\vspace{-2ex}

\hfill $ \spadesuit$

\subsection{Proof of Theorem~\ref{t6}}

The lower bound is a trivial consequence of Theorem~\ref{t5}.
Alternatively, it follows  from Theorem~\ref{t4} with $v_1=\ve $,
$v_2 = v - \ve $ and then letting $\ve \to 0$.

 To establish the complementary upper bound we fix $v>1$ and
without loss of generality  assume that $f$ satisfies
(\ref{e:011}) on $I_0$. The case that $v=1 $ is trivial. Now fix
$\ve $ such that
\begin{equation*}
0 \ < \ve \ < \ \min\{  1/(1+ v ), \, 1/5 \}    \hspace{20mm} {\rm
and } \hspace{20mm} v - \ve \ > \  1  \ \ .
\end{equation*}
The following inclusions readily follow from the definitions of
$\cSM_2(v)$ and $\cS_2(v_1,v_2)$ :
\begin{equation}\label{e:032}
\cSM_2(v) \ \subset \  \cS_2(v-\ve,0) \ \cup \  \cS_2(0,v-\ve) \
\cup \ \bigcup_{t=-t_0}^{t_0} \cS_2(v_1(t),v_2(t))  \ \subset  \
\cSM_2(v-\ve) \ ,
\end{equation}
where $t_0 > 0 $ is the unique positive integer satisfying  $
\frac{v}{2\ve} - \frac{3}{2} \leq t_0 < \frac{v}{2\ve} -
\frac{1}{2} \ $ and
\begin{equation*}
\textstyle { v_1(t) \ := \ \frac{v}{2} - \frac{(2t+1)\ve}{2}
\hspace{20mm} v_2(t) \ := \ \frac{v}{2} + \frac{(2t-1)\ve}{2} \ \
. }
\end{equation*} The required upper bound will follow on
establishing the corresponding upper bounds for the sets `between
the inclusions' of (\ref{e:032}). For this we will repeatedly
apply Theorem~\ref{t4}.

\medskip

 First, consider the sets $\cS_2(v_1(t),v_2(t))$ for $t\ge0$ (the
case $t<0$ is similar).  So, $v_1(t) \leq v_2(t)$.  Assume for the
moment that $v_1(t) < 1$.  Then $v_1(t) \in (0,1)$ and since $
v_1(t) + v_2(t) = v - \ve > 1 $ we have via  Theorem \ref{t4} that
$$
\dim  \, \cS_2(v_1(t),v_2(t))  \, \cap \, \cC_f \ =  \
\frac{2-v_1(t)}{1+v_2(t)} \ = \ \frac{2-v_1(t)}{1+v - \ve -
v_1(t)}\ \leq \ \frac{2}{1+v - \ve } \ \ .
 $$
Now suppose that $v_1(t)  \geq 1$.  It follows from the definition
of $v_1(t)$ that $v > 2 $. Trivially, $\cS_2(v_1(t),v_2(t)) \subset
\cS_2(1-\ve,v/2) $. Now, $v/2 >  1-\ve >0 $ and on applying Theorem
\ref{t4} we have that
$$
\dim \, \cS_2(v_1(t),v_2(t)) \, \cap \, \cC_f \ \leq  \ \dim \,
\cS_2(1-\ve,v/2) \, \cap  \, \cC_f \ =  \ \frac{2+2\ve}{2+v} \ \le
\ \frac{2}{1+v} \ .
$$
Next, we consider the set $ \cS_2(v-\ve,0)  $  -- the case of
$\cS_2(0,v-\ve)$ is similar.   By definition,
 $\cS_2(v-\ve,0)=\cS_1(v-\ve)\times\R$ and so
$$
\dim \, \cS_2(v-\ve,0) \, \cap \, \cC_f \ \le \ \dim\cS_1(v-\ve) \
= \ \frac{2}{1+v-\ve} \  .
$$

The upshot is that
$$
\dim \,  \cSM_2(v) \,  \cap \, \cC_f \  \leq  \ \max\left\{
\frac{2}{1+v -\ve} , \frac{2}{1+v} \right\}  \  =  \ \frac{2}{1+v
-\ve} \ ,
$$

\noindent and since $\ve$  can be made arbitrarily small the
required upper bound follows.

\vspace{-2ex}

\hfill $ \spadesuit$

\section{Final remarks: the dual form of approximation}

In view of Khintchine's transference principle \cite{Spr79},
Theorem KM can be reformulated for the dual form of approximation:

\begin{theoremKMd}\label{KM98d}
Let $\cM$ be a non-degenerate manifold in $\R^n$. Then for any
$v>1$ for almost every point $(y_1,\dots,y_n)\in \cM$ the
inequality
\begin{equation}\label{e:033}
 \|a_1y_1+\ldots+a_ny_n\| < \Pi_+(\vv a)^{-v}
\end{equation}
has only finite number of solutions $\vv
a=(a_1,\dots,a_n)\in\Z^n$, where
$$
\Pi_+(\vv a):=\prod_{i=1}^n\max\{1,|a_i|\}\,.
$$
\end{theoremKMd}

The  problems S1 and S2 considered in \S\ref{simonman}  above can
therefore be reformulated for the dual form of approximation.
Given an approximating   function $\psi$,  consider the inequality
\begin{equation}\label{e:034}
 \|a_1y_1+\ldots+a_ny_n\| < \psi(\Pi_+(\vv a))\ .
\end{equation}
Let
$$
\cLM_n(\psi):=\{\vv y\in\R^n:\text{ (\ref{e:034}) holds for
infinitely many }\vv a\in\Z^n\}
$$
and
$$ \cLM_n(v):=\cSM(q\mapsto q^{-v}):=\{\vv y\in\R^n:\text{ (\ref{e:033}) holds for
infinitely many }\vv a\in\Z^n\}\,.
$$

\noindent{\bf Problem D1\,:} Given a non-degenerate manifold
$\cM\subset\R^n$ and $v>1$, what is the Hausdorff dimension of
$\cLM(v)\cap \cM$ ?

Note that above  theorem of Kleinbock and Margulis only implies
that the Lebesgue measure of $\cLM(v)\cap \cM$ is zero.

\noindent{\bf Problem D2\,:} Given a non-degenerate manifold
$\cM\subset\R^n$ and an approximating  function $\psi$, what is
the weakest condition under which $\cLM(\psi)\cap \cM$ is of
Lebesgue measure  zero ?

\vspace{2ex}

 Regarding Problem D1 the following general lower
bound can be established:

\begin{theorem}\label{t7}
Let $\cM$ be arbitrary manifold in $\R^n$. Then for any $v>1$
\begin{equation}\label{e:035}
\dim \cM\cap\cLM_n(v) \ \ge \  \dim \cM-1+\dfrac{2}{1+v}\,.
\end{equation}
\end{theorem}

The proof of Theorem~\ref{t7} follows the same line of reasoning as
that of Theorem~\ref{t5} and is left to the reader. It is highly
likely that the inequality given by  (\ref{e:035}) is in fact an
equality. For $n=2$, that this is indeed the case is easily verified
by modifying the arguments of \cite{Yu81}. However, the general case
($n \geq 3$) seems to be a difficult problem.

%

\vspace{2ex}

 Regarding Problem D2 a general Khintchine-Groshev
type theorem for convergence has been established in \cite{BKM01}.
This states that  in Theorem \ref{KM98d} above one can  replace
(\ref{e:033}) with (\ref{e:034}) whenever
$$
\sum_{h=1}^\infty\psi(h)\log^{n-1}h \ < \ \infty \ .
$$

The divergence counterpart  remains an open problem even for
planar curves.

\vspace{4ex}

\noindent{\bf Acknowledgements.}  SV would like to thank  Ayesha
(Dorothy) and Iona (Tinman)  for allowing him into their wonderful
world of Oz and for constantly restuffing his  straw.   Also, he'd
like to thank Bridget (Wicked Witch) for her support and
friendship. Finally, many thanks to Geraldine and Peter for their
generosity and  introducing us creatures from Oz  to Rasquera.

\vspace{5mm}

\noindent Victor V. Beresnevich:  Institute of Mathematics, Academy
of Sciences of Belarus,

\vspace{-2mm}

  ~ \hspace{30mm}
220072, Surganova 11, Minsk, Belarus.


 ~ \hspace{30mm} e-mail: beresnevich@im.bas-net.by

\vspace{5mm}

\noindent Sanju L. Velani: Department of Mathematics, University of
York,

\vspace{-2mm}

 ~ \hspace{19mm}  Heslington, York, YO10 5DD, England.


 ~ \hspace{19mm} e-mail: slv3@york.ac.uk

\end{document}